\documentclass[11pt,a4paper]{article}
\usepackage[margin=1in]{geometry}
\usepackage{latexsym, amsfonts, amsmath}
\newcommand{\be}{\begin{equation}}
\newcommand{\ee}{\end{equation}}
\newcommand{\bea}{\begin{eqnarray}}
\newcommand{\eea}{\end{eqnarray}}
\newcommand{\bean}{\begin{eqnarray*}} 
\newcommand{\eean}{\end{eqnarray*}}

\newcommand{\brray}{\begin{array}}
\newcommand{\erray}{\end{array}}
\newcommand{\ben}{\begin{equation}{nonumber}}
\newcommand{\een}{\end{equation}{nonumber}}

\newtheorem{dfn}{Definition}[section]
\newtheorem{thm}[dfn]{Theorem}
\newtheorem{lmma}[dfn]{Lemma}
\newtheorem{ppsn}[dfn]{Proposition}
\newtheorem{crlre}[dfn]{Corollary}
\newtheorem{xmpl}[dfn]{Example}
\newtheorem{rmrk}[dfn]{Remark}
\newcommand{\bdfn}{\begin{dfn}}
\newcommand{\bthm}{\begin{thm}}
\newcommand{\blmma}{\begin{lmma}}
\newcommand{\bppsn}{\begin{ppsn}}
\newcommand{\bcrlre}{\begin{crlre}}
\newcommand{\bxmpl}{\begin{xmpl}}
\newcommand{\brmrk}{\begin{rmrk}}
\newcommand{\edfn}{\end{dfn}}
\newcommand{\ethm}{\end{thm}}
\newcommand{\elmma}{\end{lmma}}
\newcommand{\eppsn}{\end{ppsn}}
\newcommand{\ecrlre}{\end{crlre}}
\newcommand{\exmpl}{\end{xmpl}}
\newcommand{\ermrk}{\end{rmrk}}

\newcommand{\IC}{\mathbb{C}}

\newcommand{\IT}{\mathbb{T}}

\newcommand{\cla}{{\cal A}}
\newcommand{\clb}{{\cal B}}
\newcommand{\clc}{{\cal C}}

\newcommand{\clh}{{\cal H}}
\newcommand{\cli}{{\cal I}}

\newcommand{\cll}{{\cal L}}

\newcommand{\clq}{{\cal Q}}

\def\a*{{\cal A}_{h,*}}
\def\B{{\cal B}(h)}
\def\B1{{\cal B}_1(h)}
\def\b{{\cal B}^{\rm s.a.}(h)}
\def\b1{{\cal B}^{\rm s.a.}_1(h)}

\newcommand{\ot}{\otimes}

\newcommand{\raro}{\rightarrow}

\def \qed {$\Box$}

\def\a*{{\cal A}_{h,*}}
\def\B{{\cal B}(h)}
\def\B1{{\cal B}_1(h)}
\def\b{{\cal B}^{\rm s.a.}(h)}
\def\b1{{\cal B}^{\rm s.a.}_1(h)}

\begin{document}
\begin{center}
{\Large{\bf A note on the injectivity of action by compact quantum groups on a class of  $C^{\ast}$-algebras}}\\ 
{\large {\bf Debashish Goswami\footnote{Partially supported by J C Bose Fellowship from D.S.T. (Govt. of India) and also acknowledges
 the Fields Institute, Toronto for providing hospitality for a brief stay when a small part of this work was done.} and \bf Soumalya Joardar \footnote{Acknowledges support from CSIR.}}}\\
Indian Statistical Institute\\
203, B. T. Road, Kolkata 700108\\
Email: goswamid@isical.ac.in,\\
Phone: 0091 33 25753420, Fax: 0091 33 25773071\vspace{2mm}\\
{\bf {\large Dedicated to Paul F. Baum on the occasion of his eightieth birthday.}}\\
\end{center}
\begin{abstract}
We give some sufficient conditions for the injectivity of actions of compact quantum groups on $C^{\ast}$-algebra. As an application, 
 we prove that any faithful smooth  action by a compact quantum group on a compact smooth (not necessarily connected) manifold is 
 injective. A similar result is proved for actions on $C^{\ast}$- algebras obtained by Rieffel-deformation of compact, smooth manifolds. 
\end{abstract}
 \section{Introduction}
 Quantum groups are natural generalization of groups and they are used as `generalized symmetry objects' in mathematics and physics. 
  The pioneering work by  Drinfeld
and Jimbo (\cite{drinfeld}, \cite{drinfeldICM}, \cite{jimbo1}, \cite{jimbo}) and others  gave the formulation of  quantum groups 
 in the algebraic setting as Hopf algebras typically obtained by 
 deformations of the universal enveloping algebras of semisimple Lie algebras. This led  to a deep and successful theory having connections with
physics, knot theory, number theory, representation theory etc. On the other hand,   S.L. Woronowicz (see, e.g. \cite{Pseudogroup})  approached 
 it from a point of view of harmonic analysis on locally compact groups and came up with 
 a set of axioms for defining compact quantum groups (CQG for short) as
the generalization of compact topological groups. In this note, we will restrict ourselves to the framework of compact quantum groups only. 
 It is natural to define quantum analogue of group action on spaces. This can be done in different ways: for the purely algebraic approach, this 
  is defined as a co-action  of Hopf algebra. In the analytic theory, there are $C^{\ast}$ and von Neumann algebraic notions of action. We'll be concerned 
   with Podles' formulation of $C^{\ast}$-action of compact quantum groups on $C^{\ast}$- algebras. A subtle point about this definition is that it does not assume 
    the injectivity of the action. We mention here that some authors (e.g. \cite{free}) indeed prefer to include injectivity in the definition 
     of a $C^{\ast}$-action, but this is not a universal practice. In fact, in \cite{Act_com}, Soltan discussed examples of non-injective $C^{\ast}$-actions. On the other hand,  group actions on 
     spaces are always injective. Injectivity also follows in the algebraic setting for Hopf algebra co-actions as well as for von Neumann algebraic notion of 
      actions of (von Neumann algebraic) quantum groups. Thus, it is an interesting and important problem to give sufficient conditions for injectivity of 
       $C^{\ast}$-action of compact quantum groups in the sense of Podles. This is the aim of this short note.  We'll consider CQG actions on $C(M)$ and their Rieffel deformations where 
        $M$ is a compact smooth manifold. Under a smoothness condition on the action on $C(M)$ ( in the sense of  \cite{gafa})  we can prove injectivity. 
        For  the Rieffel-deformation of classical manifolds, we prove injectivity under a natural analogue of smoothness and compatibility of the action with the canonical toral action.

\section{ Preliminaries}
In this paper all the
Hilbert spaces are over $\mathbb{C}$ unless mentioned otherwise. For a vector space $V$, $V^{'}$ stands for its algebraic dual.
$\oplus$ and $\ot_{\rm alg}$ will denote the algebraic direct sum and
algebraic tensor
product respectively. On the other hand, the minimal $C^{\ast}$- algebra tensor product and tensor product of Hilbert spaces as well as Hilbert modules will be denoted by $\ot$. In particular, we consider Hilbert modules of the form $\clh \ot \clc$ where $\clc$ is a $C^*$ algebra. 
We shall denote the $C^{\ast}$-
algebra of bounded operators on a Hilbert space $\clh$ by $\clb(\clh)$ and the
$C^{\ast}$- algebra
of compact operators on $\clh$ by $\clb_{0}(\clh)$. $Sp$, $\overline{Sp}$ 
stand for the
linear span and closed linear span of elements of a vector space respectively, whereas ${\rm Im}(A)$ denotes the image of a linear map $A$. Given a group action $\gamma$ on a locally convex space $Z$ we denote the fixed point subspace by $Z^\gamma$.

  We call a
locally convex space Fr\'echet if the family of seminorms is countable
(hence the space
is metrizable) and the space is complete with respect to the metric given by the family
of seminorms. There are many ways to equip the algebraic
tensor product of two locally convex spaces with a locally convex topology. Let
$E_{1}$, $E_{2}$ be two
locally convex spaces with family of seminorms $\{||.||_{1,i}\}$ and
$\{||.||_{2,j}\}$ respectively. Then one wants a family $\{||.||_{i,j}\}$ of
seminorms for $E_{1}\ot_{\rm alg} E_{2}$ such that $||e_{1}\ot
e_{2}||_{i,j}=||e_{1}||_{1,i}||e_{2}||_{2,j}$. The problem is that such a choice
is far from unique and there is a maximal and a minimal choice giving the
projective and injective tensor product respectively. Let us denote the projective tensor product by $E_1 \hat{\ot} E_2$. A Fr\'echet locally convex space is called nuclear if its projective
and
injective tensor products with any other Fr\'echet space coincide as a
locally convex space. It is known that closed subspaces and quotients by closed subspaces of a nuclear Fr\'echet space are again nuclear. 
 We do not go into further details
of this topic
here but refer the reader to \cite{TVS} for a comprehensive discussion. 
Furthermore if the space is a $\ast$ algebra then we demand that its $\ast$
algebraic structure is compatible with its locally convex topology i.e.
the involution $\ast$ is continuous and multiplication is jointly continuous with
respect to the topology. Projective and
injective tensor product of two such topological $\ast$ algebras are again
topological $\ast$ algebra. We shall mostly use unital $\ast$
algebras. Henceforth all the topological $\ast$-algebras will be unital unless
otherwise mentioned. Consider a locally convex algebra $\cla$ for which each of  the defining  seminorms, say $\| \cdot \|_i$, satisfies 
  \be \label{1} \| xy\|_i \leq C_i\| x \|_i \| y\|_i\ee for some constant $C_i$ and all $x,y \in \cla$. Then it is easy to see from the definition of 
   the projective tensor product that the algebra multiplication map (say $m$) lifts to a continuous map from $\cla \hat{\ot} \cla$ to $\cla$.

We mainly need a particular class of nuclear locally convex $\ast$-algebra, which is $C^\infty(M)$, where $M$ is any compact smooth manifold.  The natural Fr\'echet topology on $C^{\infty}(M)$ is 
given by the seminorms of the form $p^{U, K,\alpha}$,
\begin{displaymath}
p^{U,K,\alpha}(f)={\rm sup}_{x\in K}|\partial^{\alpha}(f)(x)|,
\end{displaymath}
where $K$ is a compact subset contained in the domain of some coordinate chart $(U,(x_{1},...,x_{n}))$, $\alpha=(i_{1},...,i_{k})$ a multi index and 
$\partial^{\alpha}=\frac{\partial}{\partial x_{i_{1}}}...\frac{\partial}{\partial x_{i_{k}}}$, $i_j \in \{1, \ldots, n\}$.
We can similarly define a Fr\'echet topology on $C^{\infty}(M,E)$, the space of smooth $E$-valued functions on $M$ for any Fr\'echet space $E$. We refer the reader to \cite{gafa} for more details.
One can verify condition (\ref{1}) for the family of seminorms defining the Fr\'echet algebra  $C^\infty(M,\cla)$ where $\cla$ is a Banach algebra, by Leibniz rule. 
   
\subsection{Compact quantum groups and $C^{\ast}$-actions}
\bdfn
A compact quantum group (CQG for short) is a  unital $C^{\ast}$-algebra $\clq$ with a
coassociative coproduct 
(see \cite{VanDaele}, \cite{Pseudogroup}) $\Delta$ from $\clq$ to $\clq \ot \clq$  
such that each of the linear spans of $\Delta(\clq)(\clq\ot 1)$ and that
of $\Delta(\clq)(1\ot \clq)$ is norm-dense in $\clq\ot \clq$.
\edfn

 A unitary representation of a CQG $(\clq,\Delta)$ on a Hilbert
space $\clh$ is a unitary $U \in \cll(\clh \ot \clq)$ such that the $\mathbb{C}$-linear map $V$ from $\clh$ to the Hilbert module
$\clh \ot
\clq$ given by $V(\xi)=U(\xi \ot 1)$ satisfies 
 $(V \ot {\rm id})\circ V=({\rm id} \ot \Delta)\circ V.$
Here, 
the map $(V \ot {\rm id})$ denotes the extension of $V \ot {\rm id}$ to the completed tensor product 
  $\clh \ot \clq$ which exists as $V$ is an isometry.

Every  CQG $\clq$ contains   a canonical dense unital $\ast$-subalgebra
$\clq_0$ of $\clq$ on which linear maps $\kappa$ and 
$\epsilon$ (called the antipode and the counit respectively) are defined making the above subalgebra a Hopf $\ast$ algebra. In fact, this is  the algebra generated by the `matrix coefficients' of
the (finite dimensional) irreducible non degenerate representations (see \cite{VanDaele} )
 of the CQG. The antipode is an anti-homomorphism and also satisfies $\kappa(a^{\ast})=(\kappa^{-1}(a))^{\ast}$ for $a \in \clq_{0}$.
 
 It is known that there is a unique state $h$ on a CQG $\clq$ (called the Haar
state) which is bi invariant in the sense that $({\rm id} \ot h)\circ
\Delta(a)=(h \ot {\rm id}) \circ \Delta(a)=h(a)1$ for all $a$. The Haar state
need not be faithful in general, though it is always faithful on $\clq_0$ at
least. Given the Hopf $\ast$-algebra $\clq_{0}$, there can be several CQG's
which have this $\ast$-algebra as the Hopf $\ast$-algebra generated by the
matrix elements of finite dimensional representations. We need two of such CQG's: the reduced and the universal one. By definition, the reduced 
 CQG $\clq_r$ is the image of $\clq$ in the GNS representation of $h$, i.e. $\clq_r=\pi_r(\clq)$, $\pi_r: \clq \raro \clb(L^2(h))$ is the GNS representation. 
 
 There also exists a
largest such CQG $\clq^{u}$, called the universal CQG corresponding to
$\clq_{0}$. It is obtained as the universal enveloping $C^{\ast}$-algebra of
$\clq_{0}$. We also say that a CQG $\clq$ is universal if $\clq=\clq^{u}$. Given two CQG's $(\clq_{1},\Delta_{1})$ and $(\clq_{2},\Delta_{2})$, a $\ast$ homomorphism $\pi:\clq_{1}\raro\clq_{2}$ is said to be a CQG morphism if $(\pi\ot\pi)\circ\Delta_{1}=\Delta_{2}\circ\pi$ on $\clq_{1}$. In case $\pi$ is surjective, $\clq_{2}$ is said to be a quantum subgroup of $\clq_{1}$  and denoted by $\clq_{2}\leq\clq_{1}$.

\bdfn
We say that a CQG $(\clq, \Delta)$ (co)-acts on a (unital) $C^*$-algebra $\clc$ if there is a unital $\ast$-homomorphism $\alpha: \clc \raro \clc \ot \clq$ such that 
$(\alpha \ot {\rm id})\circ \alpha=({\rm id} \ot \Delta)\circ \alpha$, and the linear span of $\alpha(\clc)(1 \ot \clq)$ is norm-dense in $\clc \ot \clq$.
\edfn
 In Woronowicz theory, it is customary to drop `co', and call the above co-action simply `action' of the CQG on the $C^*$-algebra. Let us adopt this convention for the rest of the note. 
    An action  $\alpha$ on $\clc$ is called faithful if the $\ast$-subalgebra generated by $\{ (\omega \ot {\rm id})(\alpha(b)) \}$, where $b \in \clc$ and $\omega$ varying over the 
   set of bounded linear functionals on $\clc$, is dense in $\clq$.
   
   Given an action $\alpha$, we define $\alpha_r=({\rm id} \ot \pi_r)\circ \alpha$ and call it the reduced action. If the Haar state is faithful on $\clq$, we have $\alpha=\alpha_r$.
      \bdfn
 We call an action $\alpha$ of a CQG $\clq$ on a unital $C^{\ast}$-algebra $\clc$ to be implemented by a unitary representation $U$ of $\clq$ in $\clh$, say, if
  there is a faithful representation $\pi : \clc \raro \clb(\clh)$ such that $U (\pi(x) \ot 1)U^{\ast}=(\pi \ot {\rm id})(\alpha(x))$ for all $x \in \clc$.
  \edfn
  It is clear that if an action is implemented by a unitary representation then it is one-to-one. In fact, as $({\rm id} \ot \pi_r)(U)$ gives a unitary representation of $\clq_r$ in $\clh$
   and the `reduced action' $\alpha_r:=({\rm id} \ot \pi_r)\circ \alpha$ of $\clq_r$ is also implemented by a unitary representation, it follows that even 
    $\alpha_r$ is one-to-one. We see below that this is actually equivalent to implementability by unitary representation:
    \blmma
    \label{unitary_impl}
    Given an action $\alpha$ of $\clq$ on a unital separable $C^{\ast}$-algebra $\clc$, the following are equivalent:\\
    (a) There is a faithful positive functional $\phi$ on $\clc$ which is invariant w.r.t. $\alpha$, i.e. $(\phi \ot {\rm id})(\alpha(x))=\phi(x)1_\clq$ for all $x \in \clc$.\\
    (b) The action is implemented by some unitary representation.\\
    (c) The reduced action $\alpha_r$ of $\clq_r$ is injective.
\elmma
{\it Proof:}\\
If (a) holds, we consider $\clh$ to be the GNS space of the faithful positive functional $\phi$. The GNS representation $\pi$ is faithful, and the linear map $V$ defined by 
 $V(x):=\alpha(x)$ from $\clc \subset \clh=L^2(\clc,\phi)$ to $\clh \ot \clq$ is an isometry by the invariance of $\phi$. Thus $V$ extends to $\clh$ and it is easy to check that 
  it induces a unitary representation $U$, given by $U(\xi \ot q)=V(\xi)q$, which implements $\alpha$.\\
   \indent We have already argued $(b) \Rightarrow (c)$, and finally, if (c) holds, we choose any faithful state say $\tau$ on the separable $C^{\ast}$ algebra $\clc$ and take $\phi(x)=(\tau \ot h)(\alpha_r(x))$,
   which is faithful as $h$ is faithful on $\clq_r$ and $\alpha_r$ is injective. It can easily be verified that $\phi$ is $\alpha$-invariant on the dense subalgebra $\clc_0$ mentioned before, and hence on the whole of 
    $\clc$.\qed
    
   We'll also need the following facts about actions on commutative $C^{\ast}$-algebras. 
   \bppsn
   \label{abc}
   If a CQG $\clq$ acts faithfully on $C(X)$ for some compact metrizable  space $X$, then $\clq$ is separable and it is also of Kac type, i.e. $\kappa$ is norm-bounded 
    on $\clq_r$, $\kappa^2={\rm id}$ and the Haar state is tracial. 
    \eppsn
    {\it Proof}\\
    Note that $X$ is second countable and hence $C(X)$ is separable. Choose a countable dense set of points $\{ x_i,~i=1,2,\ldots \}$ and a countable norm-dense subset $\{ f_n ,~n=1,2,\ldots\}$ 
     of $C(X)$. It follows from faithfulness of the action $\alpha$ (say) that $\clq$ is generated as a $C^{\ast}$-algebra by the countable set $\{ \alpha(f_n)(x_i),~i,n=1, 2,\ldots \}$, hence it is separable.\\
     \indent For the proof of the Kac conditions, see \cite{Huichi}.\qed
 \section{Smooth actions are injective}
Let $M$ be a compact smooth manifold. Let us recall the definition of smooth CQG action on it from \cite{gafa}.
\bdfn 
We say that an action $\alpha$ of a CQG $\clq$ on $C(M)$ is smooth if 
 $\alpha$ maps $C^\infty(M)$ into $C^\infty(M, \clq)$ and   $Sp \ \alpha(C^\infty(M))(1\ot \clq)$ is dense in $C^\infty(M, \clq)$ in the Frechet topology.\edfn

 \bthm
\label{smooth_bdd_counit}
If $\clq$ has a faithful smooth action $\alpha$ on $C^\infty(M)$, where $M$ is compact manifold, then for every fixed $x \in M$ there is a well-defined, $\ast$-homomorphic 
 extension $\epsilon_x$ of the counit map
 $\epsilon$ of $\clq_0$  to the unital $\ast$-subalgebra $\clq^\infty_x:=\{ \alpha_r(f)(x):~f \in C^\infty(M)\}$ satisfying $\epsilon_x(\alpha_r(f)(x))=f(x)$, where $\alpha_r$ is 
  the reduced action discussed earlier.\ethm 
 {\it Proof:}\\
 Adapting the arguments of \cite{Podles} and \cite{free}, we can get a  the Fr\'echet dense subalgebra $\clc_0$ of $C^\infty(M)$ on which $\alpha$ restricts to an algebraic co-action of $\clq_0$. 
 For example, $\clc_0$ may be chosen as the Peter-Weyl subalgebra in the sense of \cite{free}. 
 Replacing $\clq$ by $\clq_r$ we can assume without loss of generality that $\clq$ has faithful Haar state and $\alpha=\alpha_r$. 
 In this case $\clq$ will have bounded antipode $\kappa$ (by Proposition \ref{abc}). Let $\alpha_x: C^\infty(M) \raro \clq^\infty_x$ be the map defined by $\alpha_x(f):=\alpha(f)(x)$. 
 It is clearly continuous w.r.t. the Fr\'echet topology of $C^\infty(M)$ and hence the kernel  $\cli_x$ (say)
  is a closed ideal, so that the quotient, which is isomorphic to $\clq^\infty_x$, is a nuclear space. Let us consider $\clq^\infty_x$ with this topology and then by nuclearity,
   the projective and injective tensor products with $\clq$ (viewed as a
 Banach space, which is separable by Proposition \ref{abc}) coincide.  The
multiplication map $m: \clq^\infty_x \ot_{\rm alg} \clq \raro \clq$ extends to a 
continuous map (to be denoted by $m$ again) on $\clq_x^\infty \hat{\ot} \clq$. 
          It follows from the relation $({\rm id}\ot\Delta) \circ \alpha=(\alpha \ot {\rm id})\circ \alpha$ that $\Delta(\alpha_x(f))=(\alpha_x \ot {\rm id})(\alpha(f))$, i.e. $\Delta$ 
           maps $\clq_x^\infty$ to the 
           $(\alpha_x \ot {\rm id}_\clq)(C^\infty(M) \hat{\ot} \clq) \subseteq \clq_x^\infty \hat{\ot} \clq$.  Thus, the composite map 
 $\beta:=m \circ ({\rm id} \ot \kappa) \circ \Delta: \clq^\infty_x \raro \clq$ is continuous.  Clearly, this map  coincides with 
 $\epsilon(\cdot) 1_\clq$ on the Fr\'echet-dense subalgebra of $\clq^\infty_x$ spanned by elements of the form $\alpha(f)(x)$, 
  with $f$ varying in the Fr\'echet-dense  subalgebra $\clc_0$  of $C^\infty(M)$. By continuity of $\beta$, it follows that the range of $\beta$ is $\IC 1_\clq$. 
   Hence we can define $\epsilon_x$ by setting $\epsilon_x(\cdot)1_\clq=\beta(\cdot).$ This completes the proof of the lemma.\qed
 
      \bcrlre
  \label{automatic_inj}
  For any smooth action $\alpha$ on $C^\infty(M)$, the conditions of Theorem \ref{unitary_impl} are satisfied, hence the reduced  action is injective on $C(M)$.
  \ecrlre
  {\it Proof:}\\
  Replacing $\clq$ by the Woronowicz subalgebra generated by $\{ \alpha(f)(x), f \in C(M), x \in M\}$ we may assume that 
   $\alpha$ is faithful. If $\alpha_r(f)=0$ for $f \in C^\infty(M)$  then by Lemma \ref{smooth_bdd_counit} applying the extended $\epsilon$ we conclude 
    $f=0$. Now, consider any positive Borel measure $\mu$ of full support on $M$, with $\phi_\mu$ being the positive functional obtained by integration w.r.t $\mu$. 
    Let $\psi:=(\phi_\mu \ot h) \circ \alpha_r$ be the positive functional which is clearly $\alpha_r$-invariant and faithful on $C^\infty(M)$, 
    i.e. $\psi(f)=0, f \in C^\infty(M)$ and $f$ nonnegative  implies $f=0$. But  by Riesz Representation Theorem there is a positive Borel measure
     $\nu$ such that $\psi(f)=\int_M f d\nu$. We claim that $\nu$ has full support, hence $\psi$ is faithful also on $C(M)$. 
     Indeed, for any nonempty open subset $U$ of $M$ there is a nonzero positive $f \in C^\infty(M)$, with $0 \leq f \leq 1,$ and support of $f$ is contained in $ U$.   
     By faithfulness of $\psi$ on $C^\infty(M)$ we get $0 < \psi(f) =\int_U f d\nu \leq \nu(U).$\qed
 
\brmrk
In a recent work \cite{final}, Goswami has proved that any CQG which admits a faithful smooth action on a compact connected smooth manifold must be isomorphic to $C(G)$ for some 
 group $G$ and the CQG action becomes $G$-action. Hence it is in particular injective. However, the result of the present note is applicable to a  possibly disconnected manifold. Moreover, 
  the proof of injectivity given here is rather short and direct. 
\ermrk

\section{Smooth action on Rieffel deformation} 
Let us now consider CQG actions on noncommutative $C^{\ast}$-algebras. 
Rieffel deformation (see \cite{rieffel}) is a well-known and very useful procedure to obtain interesting noncommutative $C^{\ast}$-algebras from the commutative ones. In particular, given  a smooth 
 compact Riemannian manifold $M$ equipped with a toral subgroup $T \cong \IT^n \subseteq ISO(M)$, one can construct a family of (typically noncommutative) $C^{\ast}$-algebras 
  $C(M)_\theta$ indexed by $n \times n$ skew symmetric matrices $\theta$. There is a similar procedure (see \cite{wang_def}) for deforming a CQG 
    with some toral quantum subgroup of rank $n$ inducing an action of torus of rank $2n$ combining left and right action by the elements of the $n$-toral subgroup. 
     In this case, one gets a CQG by retaining the same co-algebra structure as the original one but changing  the algebra structure. This will be called the Rieffel-Wang deformation of $G$.  
     
     Let $\cla_\theta$ be the noncommutative $n$-torus, which is the universal $C^{\ast}$-algebra 
   generated by unitaries $U_1, \ldots U_n$ satisfying the commutation relations $U_j U_k={\rm exp}(2\pi i \theta_{jk})U_k U_j$, where $\theta=(( \theta_{jk} )).$ 
    Given a unital $C^{\ast}$-algebra $\clc$ with a $\IT^n$-action given by $\beta_z$ (say), the deformed $C^*$-algebra $\clc_\theta$ can be described in two alternative ways: either in the original picture of Rieffel where one defines a new, twisted multiplication on 
     the spectral algebra for the toral action and then considers appropriate $C^{\ast}$-completion, or as   in \cite{connes_dubois}, identifying 
  $\clc_\theta$  with the fixed point subalgebra $(\clc \ot \cla_\theta)^{\beta \ot v^{-1}}$  where    
  $v_z$ denotes the canonical toral action on $\cla_\theta$ satisfying $v_z(U_i)=z_iU_i$ for all $i$. We have a `dual' $T$-action on $\clc_\theta$ which is the restriction of $({\rm id} \ot v)$ on $\clc \ot \cla_\theta$. 
  
    Given a CQG $\clq$ and a quantum subgroup of $\clq$ isomorphic with $T=\IT^n$, with the corresponding 
       surjective CQG morphism $\pi : \clq \raro C(T)$, we can define a left and a right $\IT^n$-action, say $\chi^l_z,\chi^r_z$ respectively, by setting 
       $\chi^l_z=({\rm id} \ot ({\rm ev}_z \circ \pi))\circ \Delta$ and $\chi^r_z=(( {\rm ev}_z \circ \pi) \ot {\rm id})\circ \Delta.$ Using this, we have a $\IT^{2n}$-action 
       $\chi_{z,w}=\chi^l_{z}\chi^r_{w}$ on $\clq$ and the corresponding deformed CQG   is the $C^{\ast}$-algebra $\clq_{\theta \oplus (-\theta)}.$

        We have the following from Theorem 3.11 of \cite{bh_gos}.
        \blmma
        \label{action}
        Let $\clc$ be a unital $C^{\ast}$-algebra equipped with a $\IT^n$-action given by $\ast$-automorphism $\beta_z$,  $\clq$ be a reduced CQG with an action $\alpha$ of $\clq$ on $\clc$ and a quantum subgroup of $\clq$ isomorphic with $\IT^n$ as above (with the corresponding 
        morphism $\pi$) satisfying $\beta_z:=({\rm id}\ot ({\rm ev}_z \circ \pi))\circ \alpha$.  Then we have an action $\alpha_\theta$ of $\clq_{\tilde{\theta}}$ ($\tilde{\theta}=
        \theta \oplus (-\theta)$) on $\clc_\theta$. Here the deformation of $\clc$ is taken w.r.t. the action $\beta$. 
        \elmma

   Consider now $\clc=C(M)$, $M$ being a compact smooth Riemannian manifold equipped with $\IT^n$ action, which also induces a $\IT^n$ action (say $\beta$) on $C(M)$. 
    Let $\gamma=\beta \ot v^{-1}$ as before and let 
     us call the subalgebra $C^\infty(M, \cla_\theta)^\gamma  \subset C(M, \cla_\theta)^\gamma\equiv 
      (C(M) \ot \cla_\theta)^\gamma=C(M)_\theta$ the `smooth subalgebra' and call an action $\alpha$ on $C(M)_\theta$ by a CQG $\clq$ 
    to be 
     smooth if it maps the above smooth subalgebra into $C^\infty(M, \cla_\theta \ot \clq)^{(\gamma \ot {\rm id})}$ and the linear span of  $\alpha(F)(1 \ot \clq)$, $F \in C^\infty(M,\cla_\theta)^\gamma$ is 
      dense in $C^\infty(M,\cla_\theta \ot \clq)^{(\gamma \ot {\rm id})}$. Now, we can state and prove the following:
      \bthm
      Let $M$ be  as above and let  $\alpha$ be a  smooth action of a CQG on $C(M)_\theta$ in the above sense. Moreover, assume that there is a quantum subgroup of $\clq$ isomorphic with $T=\IT^n$, given by a 
       surjective CQG morphism $\pi : \clq \raro C(T)$ such that $({\rm id}\ot ({\rm ev}_z \circ \pi))\circ \alpha$ coincides with the canonical `dual' $T$-action  on $C(M)_\theta$. Then the 
        action (in fact, the reduced one) is injective.
      \ethm
{\it Proof:}\\ We only very briefly sketch the proof.  As before, assume without loss of generality that the CQG is reduced. 
 It follows from the proof of Theorem 3.11 of \cite{bh_gos} that the action $\alpha_{-\theta}$ of $\clq_{\tilde{\theta}}$ on  $(C(M)_\theta)_{-\theta}
\cong C(M)$   is  smooth. We note that the word `smooth' in the statement of Theorem 3.11 of \cite{bh_gos} is used in a  sense weaker  than ours: it only means the 
 invariance of the smooth algebra there. However, the dense  subalgebra (e.g. Peter-Weyl subalgebra)of  $C(M)_\theta$ for 
  the action $\alpha$, on which $\clq_0$ (co)acts algebraically, can be identified as a vector space  with a dense subalgebra of 
   $C^\infty(M)$ on which the deformed action (which is the same as $\alpha$ a linear  map on this space) $\alpha_{-\theta}$ is algebraic. From this, 
    the Podles-type density condition follows, i.e. $\alpha_{-\theta}$ is smooth in our sense. Hence it is injective by Corollary \ref{automatic_inj}. Moreover, by that corollary and 
     Theorem \ref{unitary_impl}, we get a unitary representation of $\clq_{\tilde{\theta}}$ which implements $\alpha_{-\theta}$. But by the generalities of Rieffel-Wang (or, more general cocycle-twisted)
      deformation of CQG as in the Chapter 7 of \cite{book}, we conclude that $\alpha=(\alpha_{-\theta})_\theta$ is unitarily implemented too, where the corresponding Hilbert space and unitary 
       essentially remain the same. In particular, $\alpha$ is injective.\qed 

{\bf Ackowledgement:}\\ The first author would like to recall his fond memories of interaction with Prof. Paul Baum on several occasions: ICTP (Trieste), TIFR (India) and also in the conference in 
 the Fields Institute in 2016. Although there was no discussion on the particular topic of this note, the first author is grateful to him for encouragement and suggestions regarding Baum-Connes' conjecture 
  for quantum groups and other problems related to quantum symmetry on operator algebras. The authors wish him a long, healthy and creative  life.

\end{document}